\documentclass[preprint,12pt]{elsarticle}

\usepackage{amssymb}
\usepackage{amsmath}
\usepackage{amsthm}
\usepackage{tikz}
\usepackage{color}
\usepackage{amsmath}
\usepackage{graphics,url,color,epsfig}
\usepackage[colorlinks]{hyperref}
\AtBeginDocument{%
  \hypersetup{%
    linkcolor=blue,%
    citecolor=green,%
  }%
}
\usepackage{cleveref}

\crefname{equation}{}{}

\newtheorem{theorem}{Theorem}[section]
\newtheorem{lemma}[theorem]{Lemma}
\newtheorem{proposition}[theorem]{Proposition}
\newtheorem{corollary}[theorem]{Corollary}

\theoremstyle{definition}

\theoremstyle{remark}
\newtheorem{remark}[theorem]{Remark}

\numberwithin{equation}{section}

\journal{***}

\begin{document}

\begin{frontmatter}



\title{Liouville theorems for fractional parabolic equations \tnoteref{t1}}

\author[rvt1]{Wenxiong Chen}
\ead{wchen@yu.edu}

\author[rvt2]{Leyun Wu}
\ead{leyunwu@126.com}

\tnotetext[t1]{This research is partially supported by NSFC-12031012 and China Postdoctoral Science Foundation (No.2019M661472).}

\address[rvt1]{Department of Mathematical Sciences, Yeshiva University, New York, NY, 10033 USA}
\address[rvt2]{School of Mathematical Sciences, MOE-LSC,
Shanghai Jiao Tong University, Shanghai, PR China}

\begin{abstract}

 In this paper, we establish several Liouville type theorems for entire solutions to  fractional parabolic equations. We first obtain the key ingredients needed in the proof of Liouville theorems, such as narrow region principles and maximum principles for antisymmetric functions in unbounded domains, in which we remarkably weaken the usual decay condition $u \to 0$ at infinity with respect to the spacial variables to a polynomial growth on $u$ by constructing auxiliary functions.
Then we derive monotonicity for the solutions in a half space $\mathbb{R}_+^n \times \mathbb{R}$
 and obtain some new connections between the nonexistence of solutions in a half space $\mathbb{R}_+^n \times \mathbb{R}$ and
in the whole space $\mathbb{R}^{n-1} \times \mathbb{R}$ and therefore prove the corresponding Liouville type theorems.

To overcome the difficulty caused by the non-locality of the fractional Laplacian, we introduce several new ideas which will
become useful tools in investigating qualitative properties of solutions for a variety of non-local parabolic problems.
\end{abstract}


\end{frontmatter}

\section{Introduction}
In this paper, we establish Liouville theorems for the solutions to the following fractional  parabolic equations in both the whole space
\begin{eqnarray}\label{wholee}
\frac{\partial u} {\partial t}(x,t) +(-\Delta)^s u(x,t)= f(u(x,t)),\,\, (x, t) \in  \mathbb{R}^n \times \mathbb{R},
\end{eqnarray}
and the half space
\begin{eqnarray}\label{halfe}\left\{\begin{array}{ll}
 \frac{\partial u}{\partial t} (x, t)+( - \Delta )^s u(x, t)=f(u(x,t)), & (x, t) \in  \mathbb{R}_+^n \times \mathbb{R}, \\
u(x,t) =0,& (x,t )\notin \mathbb{R}_+^n \times \mathbb{R},
\end{array} \right. \end{eqnarray}
where $0<s<1$, $f$ is a $C^1$ function, the ranges for time variable $t$  are  $(-\infty, \infty)$, and in this case, the solutions are referred to as entire solutions.

For each fixed $t\in \mathbb{R},$ the fractional Laplacian on $x$ is defined by
\begin{eqnarray*}
(-\Delta)^s u (x,t)&=&C_{n, s} P.V. \int_{\mathbb{R}^n} \frac{u(x, t)-u(y,t)}{|x-y|^{n+2s}}dy\\
&=&C_{n, s}  \lim_{\epsilon \rightarrow 0^{+}} \int_{\mathbb{R}^n \backslash B_\varepsilon (x)} \frac{u(x, t)-u(y,t)}{|x-y|^{n+2s}}dy,
 \end{eqnarray*}
where $P.V.$ stands for the Cauchy principal value.
It is easy to see that    for $u \in C^{1,1}_{loc}\cap   {\cal L}_{2s } ,$ $(-\Delta)^su$ is well defined, where
 $$ {\mathcal L}_{2s}=\left\{u(\cdot, t) \in L^1_{loc} (\mathbb{R}^n) \mid \int_{\mathbb R^n} \frac{|u(x,t)|}{1+|x|^{n+2s}}dx<+\infty\right\}.$$
It is well-known that the fractional Laplacian is a nonlocal operator, and as $s \to 1,$ it approaches to the regular Laplacian $-\Delta.$

In this paper, we exhibit some new  monotonicity and Liouville type results for entire solutions of \eqref{wholee} and \eqref{halfe}.
It is well known that Liouville type theorems play crucial roles in the theory of PDEs, while monotonicity usually is a key tool to derive Liouville theorems. To this end, a number of systematic approaches have been established, such as the extension method (\cite{CS}), the method of moving planes (\cite{BCN1, BHR, CL, CLM, GNN}), the method of moving spheres (\cite{CLZ0, LZ1}), and the sliding methods (\cite{BN, WC2}).
From the heuristic point of view, Liouville-type theorems turn out to
be equivalent to universal (initial or final) blow-up or decay estimates (see \cite{G}, \cite{PQS}, \cite{X} and the references therein).

For the elliptic equations involving either local or nonlocal operators, there have been numerous articles that dedicated to the study of Liouville type theorems, such as
\cite{BM, D, DDWW,  GS1, GS2, HL, LWZ, PQS1, SZ} and so on.

Entire solutions to semi-linear parabolic equations involving the regular Laplacian
\begin{eqnarray*}
\frac{\partial u} {\partial t}(x,t) -\Delta u(x,t)= f(u(x,t)),\,\, (x, t) \in  \mathbb{R}^n \times \mathbb{R}
\end{eqnarray*}
play an important role in the dynamics of solutions to the Cauchy problem
\begin{eqnarray}\label{heate}\left\{\begin{array}{ll}
 \frac{\partial u}{\partial t} (x, t) - \Delta  u(x, t)=f(u(x,t)), & (x, t) \in  \mathbb{R}^n \times (0, +\infty), \\
u(x,0) =u_0(x),& x \in \mathbb{R}^n .
\end{array} \right. \end{eqnarray}
For example, the $\omega$-limits sets of bounded solutions to \eqref{heate} and global
attractors are comprised of entire solutions. For these reasons, entire solutions to reaction
diffusion equations have been widely studied (see \cite{H, MZ1, T} ).

For nonlinear parabolic equations involving local operators, such as $-\Delta$, there have also been a series results in this respect.

Bidaut-V\'{e}ron \cite{BV}, Merle and Zaag \cite{MZ1}, and Pol\'{a}\v{c}ik, Quittner and  Souplet \cite{PQS} proved that the only nonnegative bounded classical solution
of the parabolic problem
\begin{eqnarray}\label{wholee1}
\frac{\partial u}{\partial t} (x, t) - \Delta  u(x, t)=u^p(x,t), \,\, (x, t) \in  \mathbb{R}^n \times \mathbb{R}
\end{eqnarray}
is the trivial solution $u \equiv 0$ for $1<p<\frac{n(n+2)}{(n-1)^2}.$ Then Quittner generalized this result to  $1<p<\frac{n+2}{n-2}$ in \cite{Q}.
This is also true if $u^p$ is  replaced by $f(u)$, where $f$ is a decreasing continuous function and satisfies $f(c)=0$ iff $c=0.$

For the case of a half space, Kavian \cite{K},  and Levine and Meier \cite{LM} showed that
the only nonnegative bounded classical solution
of the parabolic problem
\begin{eqnarray}\label{halfe1}\left\{\begin{array}{ll}
 \frac{\partial u}{\partial t} (x, t) - \Delta  u(x, t)=u^p(x,t), & (x, t) \in  \mathbb{R}^n_+ \times \mathbb{R}, \\
u(x,t) =0,& (x,t) \in \{x_1=0\}\times \mathbb{R}^{n-1} \times \mathbb{R}
\end{array} \right. \end{eqnarray}
is the trivial solution $u \equiv 0$ for  $1<p\leq \frac{n+3}{n+1}.$ Then Pol\'{a}\v{c}ik, Quittner and  Souplet \cite{PQS} and Xing \cite{X} extended the
range of the exponent $p$ from $p\leq \frac{n+3}{n+1}$ to $p< \frac{n^2-1}{(n-2)^2}.$ Similar to the whole space case, the same conclusion holds if $u^p$ is  replaced by $f(u)$, where $f$ is a decreasing continuous function and satisfies $f(c)=0$ iff $c=0.$


So far as we aware, not much is known concerning the entire solutions to fractional parabolic equations. Due to the non-locality of the fractional Laplacian, many traditional approaches for local elliptic operators do not work anymore in the nonlocal setting. However these qualitative properties of solutions, in particular, the Liouville type theorems are definitely important tools  in the blow-up rate, a priori bounds, and optimal universal estimates of solutions to related initial and initial-boundary value nonlocal parabolic problems and so on. This is a motivation for the present paper.

In this paper, we prove Liouville type theorems for fractional parabolic problems  \eqref{wholee} and \eqref{halfe}  and establish some new connections between the solutions of \eqref{wholee} and \eqref{halfe}. To this end, we first introduce some new ideas to develop a fractional parabolic version of the method of moving planes, in which the key ingredients are {\em narrow region principles} and {\em maximum principles for anti-symmetric functions} as we will state in the following two theorems.

For $x \in \mathbb{R}^n,$ denote
$$
x=(x_1, x')
$$
where $x' \in \mathbb{R}^{n-1}$, and
$$
\Sigma_\lambda:=\left\{x \in \mathbb{R}^n \mid x_1<\lambda\right\}.
$$

\begin{theorem}\label{NR} (Narrow region principle)
Let $\Omega$ be a bounded or unbounded narrow region in $\Sigma_\lambda$, such that it is contained in $\{x \mid \lambda-2l<x_1<\lambda\}$ with small $l.$  Suppose that  $w_\lambda(x,t) \in ( C_{loc}^{1,1}(\Omega) \cap {\cal L}_{2s})\times C^1(\mathbb{R})$ is uniformly bounded with respect to $t$ and  lower semi-continuous in $x$ on $\bar{\Omega}$, and it satisfies
 \begin{eqnarray} \label{GC}
w_\lambda(x,t)\leq o(1) |x | ^\gamma, \mbox{ for any }  0<\gamma <2s \mbox{ as} \,\,| x | \to +\infty,
 \end{eqnarray}
 and
  \begin{eqnarray}\label{EENR}\left\{\begin{array}{ll}
 \frac{\partial w_\lambda}{\partial t} (x, t)+( - \Delta )^s w_\lambda(x, t)=c_\lambda(x, t) w_\lambda(x, t), & (x, t) \in  \Omega \times \mathbb{R}, \\
 w_\lambda(x,t)\geq 0,& (x,t )\in  (\Sigma_\lambda \backslash \Omega) \times \mathbb{R},\\
 w_\lambda(x^\lambda,t)=- w_\lambda(x,t),& (x,t )\in \Sigma_\lambda \times \mathbb{R}.
\end{array} \right. \end{eqnarray}

 If $c_\lambda(x, t)$ is bounded from above, then for sufficiently small $l,$ we have
\begin{equation}
w_\lambda(x, t)\geq 0, \,\, (x, t) \in  \Sigma_\lambda \times \mathbb{R}.
\label{A10}
\end{equation}
Furthermore, the following strong maximum principle holds:

Either
$$
w_\lambda(x, t)>0,\, (x,t) \in \Omega \times \mathbb{R},
$$
or
$$
w_\lambda(x,t) \equiv 0, \, (x,t) \in \Omega \times \mathbb{R}.
$$

\end{theorem}

In previous literature, for instance in \cite{CWNH},  to establish a narrow region principle in an unbounded domain,
one usually assumed that
$$ \lim_{|x| \to +\infty} w_\lambda(x,t) = 0.$$
Here in Theorem \ref{NR},  we remarkably weaken this condition and even allow $w_\lambda (x,t)$ to go to infinity and only assume its growth rate does not exceed $|x|^\gamma \,(0<\gamma <2s ).$
To this end, we consider
$$\bar{w}_\lambda(x, t)=\frac{e^{mt}w_\lambda(x,t)}{h(x)}$$
with
 a suitable auxiliary function $h(x)$ satisfying
$$
\lim_{|x| \to +\infty} h(x) = + \infty
$$
and
$$
\frac{(-\Delta)^s h(x)}{h(x)} \geq C >0.
$$
The main difficulty lies in how to seek such an $h(x)$.

For the Laplace operator, it is much easier to find such a function since we can choose $h(x)=\sin (a x_1)V(|x'|),$  where $V$ is a radial function satisfying an ordinary differential equation. However for the fractional Laplacian,  due to the nonlocality, it seems impossible to use an ODE method to find such a function $V$.
After careful calculations, we finally choose
$$
  h(x)=\left[\left(1-\frac{(x_1-(\lambda-l))^2}{l^2}\right)_+^s+1\right](1+\mid x '\mid^2)^{\frac{\gamma}{2}}.
  $$
This implies
$$
\lim_{|x|\to+\infty} \bar{w}_\lambda(x, t)=0.
$$
Consequently, when we use a contradiction argument to prove (\ref{A10}),
$\bar{w}_\lambda$ would be able to attain its negative minimum in the interior of $\Sigma_\lambda.$

For the regular Laplacian, the convenience is:
\begin{eqnarray}\label{com}
\Delta w_\lambda(x, t)=e^{-mt}(\Delta \bar{w}_\lambda(x, t)\cdot h(x)+2\nabla\bar{w}_\lambda(x, t)\cdot \nabla h(x)+\bar{w}_\lambda(x, t) \cdot\Delta h(x)).
\end{eqnarray}
At a minimum of $\bar{w}_\lambda(x, t),$ the middle term on the right hand side vanishes since $\nabla \bar{w}_\lambda(x,t)=0.$ This makes the analysis much easier. However, the fractional counter part of \eqref{com} is
$$
(-\Delta)^s w_\lambda(x, t)=e^{-mt}\left( (-\Delta)^s \bar{w}_\lambda(x, t)\cdot h(x)\right.$$
$$
\left.-2C \int_{\mathbb{R}^n} \frac{(\bar{w}_\lambda(x, t)-\bar{w}_\lambda(y, t))(h(x)-h(y))}{|x-y|^{n+2s}}dy+\bar{w}_\lambda(x, t) \cdot(-\Delta)^s h(x) \right)
$$
At a minimum of $\bar{w}_\lambda(x,t)$, the middle term on the right hand side (the integral) neither vanishes nor has a definite sign. This is the main difficulty. To circumvent it, we combine the first two terms together to derive a good estimation.

\begin{remark}\label{ellipticNR}  In \cite{CH}, the authors derive a narrow region principle for fractional elliptic equations in unbounded domains under the conditions that $w_\lambda(x)$ is bounded. Here by using our auxiliary function $h(x)$, we will be able to weaken this boundedness condition to a polynomial growth. More precisely, as a byproduct, we can prove the following:

\begin{proposition}
Let $\Omega$ be a bounded or unbounded narrow region in $\Sigma_\lambda$, such that it is contained in $\{x \mid \lambda-2l<x_1<\lambda\}$ with small $l.$  Suppose that  $$w_\lambda(x) \in  C_{loc}^{1,1}(\Omega) \cap {\cal L}_{2s}(\mathbb{R}^n)$$ is lower semi-continuous on $\bar{\Omega}, $ where $${\mathcal L}_{2s}(\mathbb{R}^n)=\{u \in L^1_{loc} (\mathbb{R}^n) \mid \int_{\mathbb R^n} \frac{|u(x)|}{1+|x|^{n+2s}}dx<+\infty\}.$$ Assume
 \begin{eqnarray} \label{EGC}
w_\lambda(x)\leq o(1) | x | ^\gamma, \mbox{ for any }  0<\gamma <2s \mbox{ as} \mid x \mid \to +\infty.
 \end{eqnarray}
 and
  \begin{eqnarray}\label{ENR}\left\{\begin{array}{ll}
( - \Delta )^s w_\lambda(x)+c_\lambda(x) w_\lambda(x)=0, & x \in  \Omega , \\
 w_\lambda(x)\geq 0,& x \in  \Sigma_\lambda \backslash \Omega ,\\
 w_\lambda(x^\lambda)=- w_\lambda(x),& x\in \Sigma_\lambda .
\end{array} \right. \end{eqnarray}

 If $ c_\lambda(x) $ is bounded from above, then we have
 $$
w_\lambda(x)\geq 0, \,\, x \in  \Sigma_\lambda.
$$
Furthermore, the following strong maximum principle holds:

Either
$$
w_\lambda(x)>0,\, x \in \Omega,
$$
or
$$
w_\lambda(x) \equiv 0, \, x \in \Omega.
$$
\end{proposition}
\end{remark}
\medskip

\begin{theorem}\label{MPA} (Maximum principle for antisymmetric functions)
Let $\Omega$ be a bounded or  unbounded region in $\Sigma_\lambda$, assume that the width of $\Omega$ in $x_1$ direction is bounded. Suppose that  $w_\lambda(x,t) \in ( C_{loc}^{1,1}(\Omega) \cap {\cal L}_{2s})\times C^1(\mathbb{R})$ is uniformly bounded with respect to $t$ and  lower semi-continuous in $x$ on $\bar{\Omega}$, and it satisfies
 \begin{eqnarray*}
w_\lambda(x,t)\leq o(1) | x| ^\gamma, \mbox{ for any }  0<\gamma <2s \mbox{ as}\,\, |x| \to +\infty,
 \end{eqnarray*}
 and
  \begin{eqnarray}\label{EMR}\left\{\begin{array}{ll}
 \frac{\partial w_\lambda}{\partial t} (x, t)+( - \Delta )^s w_\lambda(x, t)=c_\lambda(x, t) w_\lambda(x, t), & (x, t) \in  \Omega \times \mathbb{R}, \\
 w_\lambda(x,t)\geq 0,& (x,t )\in  (\Sigma_\lambda \backslash \Omega)  \times \mathbb{R},\\
 w_\lambda(x^\lambda,t)=- w_\lambda(x,t),& (x,t )\in \Sigma_\lambda \times \mathbb{R}.
\end{array} \right. \end{eqnarray}

 If
 \begin{eqnarray}\label{Clambda}
c_\lambda(x, t) \leq 0 \,\, \mbox{ or } \,\,c_\lambda(x, t) > 0 \,\mbox{ is small},\,\, (x, t) \in  \Omega \times \mathbb{R},
 \end{eqnarray}
then
$$
w_\lambda(x, t)\geq 0, \,\, (x, t) \in  \Sigma_\lambda \times \mathbb{R}.
$$
Furthermore, the following strong maximum principle holds:

Either
$$
w_\lambda(x, t)>0,\, (x,t) \in \Omega \times \mathbb{R},
$$
or
$$
w_\lambda(x,t) \equiv 0, \, (x,t) \in \Omega \times \mathbb{R}.
$$

\end{theorem}

\begin{remark}\label{ellipticMPA1}
In Theorem \ref{MPA}, if we replace the first equation in \eqref{EMR} by
\begin{eqnarray*}
 \frac{\partial w_\lambda}{\partial t} (x, t)+( - \Delta )^s w_\lambda(x, t)\geq 0, \, (x,t) \in \Omega \times \mathbb{R},
\end{eqnarray*}
then the same  conclusions as  Theorem \ref{MPA} still holds.
\end{remark}

Based on Theorem \ref{NR} and \ref{MPA}, we derive the strict monotonicity of solutions to the problem
\begin{eqnarray}\label{he}\left\{\begin{array}{ll}
 \frac{\partial u}{\partial t} (x, t)+( - \Delta )^s u(x, t)=f(u(x, t)), & (x, t) \in  \mathbb{R}_+^n \times \mathbb{R}, \\
u(x,t) =0,& (x,t )\notin \mathbb{R}_+^n \times \mathbb{R}.
\end{array} \right. \end{eqnarray}
We also find some new connections between the existence of solutions of  \eqref{he} and of the problem
\begin{eqnarray}\label{WS}
\frac{\partial u}{\partial t} (x, t) +(-\Delta)^s u(x, t)= f(u(x, t)), \,\, (x, t) \in \mathbb{R}^{n-1} \times \mathbb{R}.
\end{eqnarray}

\medskip

\begin{theorem}\label{MH} (Monotonicity in a half space)
Assume $f:[0, \infty) \to \mathbb{R}$ is a $C^1$ function with $f(0)=0, \, f'(0) \leq 0,$ and suppose
$$u \in( C_{loc}^{1,1}(\mathbb{R}_+^n) \cap C( \overline {\mathbb{R}_+^n}) \cap {\cal L}_{2s})\times C^1(\mathbb{R}). $$
Then

(i) If $u$ is a positive solution of \eqref{he} satisfying
$$u(x,t)\leq o(1) | x| ^\gamma, \mbox{ for any }  0<\gamma <2s \mbox{ as}\,\, |x| \to +\infty,$$
 then it is increasing in $x_1$ and
$$
\frac{\partial u}{\partial x_1} (x, t ) >0,\,\, (x, t) \in \mathbb{R}_+^n \times \mathbb{R}.
$$

(ii) If there is a positive bounded solution of \eqref{he}, then there exists a positive bounded solution of
\eqref{WS}.

\end{theorem}

Combining Theorem \ref{MH} (ii) with the nonexistence result in the whole space obtained in \cite{GK}, we derive the following nonexistence result in the half space.

\begin{theorem}\label{NH} (Nonexistence)
Assume $n \geq 3$, $1 < p\leq \frac{n-1+2s}{n-1}.$ Then the problem
\begin{eqnarray}\label{nhe}\left\{\begin{array}{ll}
 \frac{\partial u}{\partial t} (x, t)+( - \Delta )^s u(x, t)=u^p(x, t), & (x, t) \in  \mathbb{R}_+^n \times \mathbb{R}, \\
u(x,t) =0,& (x,t )\notin \mathbb{R}_+^n \times \mathbb{R}.
\end{array} \right. \end{eqnarray}
possesses   no nontrivial nonnegative bounded solution that is in
 $$( C_{loc}^{1,1}(\mathbb{R}_+^n) \cap C(\overline{\mathbb{R}_+^n})\cap {\cal L}_{2s})\times C^1(\mathbb{R}). $$
\end{theorem}
\medskip

Our approach can also be applied to derive Liouville type theorems for other fractional parabolic problems such as

\begin{theorem}\label{wid}
Assume $f:[0, \infty) \to \mathbb{R}$ is a $C^1$ function satisfying $f'\leq 0$,  $u \in( C_{loc}^{1,1} \cap {\cal L}_{2s})\times C^1(\mathbb{R}) $ is a nonnegative bounded solution of
\begin{eqnarray}\label{wide}
 \frac{\partial u}{\partial t} (x, t)+( - \Delta )^s u(x, t)=f(u(x, t)), & (x, t) \in  \mathbb{R}^n \times \mathbb{R},
\end{eqnarray}
Then

(i) $u$ is independent of spatial variables $x_i \,(1\leq i \leq n)$, i.e.
$$
u_t =f(u(t)), \,\, \forall \, t \in \mathbb{R}.
$$

(ii) Assume  $f(c)=0$ iff $c=0$, then the only nonnegative bounded solution of \eqref{wide} is the trivial solution $u \equiv 0.$

(iii) Assume $f(c)<0$ as $c \geq 0$, then the nonnegative bounded solution of \eqref{wide} does not exist.
\end{theorem}

Such examples of $f(u)$ are $1-e^u, \, -u^p, -e^u ...$.
\medskip

These kinds of results are also
indispensable tools in the study of blow-up rate for indefinite problems (see \cite{X}).
\medskip

For other results concerning nonlocal parabolic equations, please see \cite{JW}, \cite{GK}, \cite{MSWZZ}, and \cite{CWZ}.
\medskip

In Section 2, we derive {\em narrow region principles} and {\em maximum principles for anti-symmetric functions}. In Section 3,
we establish  {\em the monotonicity and nonexistence in half spaces}. In Section 4, we prove Theorem \ref{wid}.

\section{Maximum principles}

Before giving  the proof of narrow region principle, we first compute $\frac{(-\Delta)^sh(x)}{h(x)}.$

\begin{lemma}\label{mlem1}
Let $x=(x_1, x'), $ and
  $$
  h(x)=f(x_1)g_\gamma(x'), \,\, 0<\gamma <2s,
  $$
where
$$
f(x_1)=\left(1-\frac{x_1^2}{l^2}\right)_+^s+1,
$$
and
$$
g_\gamma(x')=(1+| x '|^2)^{\frac{\gamma}{2}}, \, \, 0<\gamma <2s.
$$
If $l>0$ is sufficiently small,  then
\begin{eqnarray}\label{eigen}
 \frac{(-\Delta)^sh(x)}{h(x)} \geq \frac{C_1}{l^{2s}}, \,\, | x_1 | <l.
\end{eqnarray}
\end{lemma}
\textbf{Proof.}
Firstly, if $| x_1 | <l,$ we have
\begin{eqnarray}\label{vf}
 (-\Delta)^s f(x_1)&=&(-\Delta)^s\left[\left(1-\frac{x_1^2}{l^2}\right)_+^s+1\right]\nonumber\\
 &=&(-\Delta)^s\left(1-\frac{x_1^2}{l^2}\right)_+^s\nonumber\\
& =&\frac{1}{l^2}(-\Delta)^s\left(1-{x_1^2}\right)_+^s\nonumber\\
&=&\frac{C}{l^2},
\end{eqnarray}
and $ (-\Delta)^s g_\gamma(x')$ is bounded. Indeed, if $| x' |=0,$ then
\begin{eqnarray*}
(-\Delta)^s g_\gamma(x')&=&C_{n,s} P.V. \int_{\mathbb{R}^n}\frac{1-(1+| y' |^2)^{\frac{\gamma}{2}}}{(| y'|^2+y_1^2)^{\frac{n+2s}{2}}} dy'dy_1\\
 &=&C_{n,s} P.V. \int_{\mathbb{R}^n}\frac{1-(1+| y'|^2)^{\frac{\gamma}{2}}}{(| y' |^2+| y' |^2 t^2)^{\frac{n+2s}{2}}} | y'| dy'dt\\
  &=&C_{n,s} P.V. \int_{\mathbb{R}^n}\frac{1-(1+| y' |^2)^{\frac{\gamma}{2}}}{|y'|^{n-1+2s}(1+ t^2)^{\frac{n+2s}{2}}} dy'dt\\
  &=&C P.V. \int_{\mathbb{R}^{n-1}}\frac{1-(1+| y'|^2)^{\frac{\gamma}{2}}}{|y'|^{n-1+2s}} dy'\\
  &=& C \int_{0}^{+\infty}\frac{1-(1+r^2)^{\frac{\gamma}{2}}}{r^{1+2s}}dr
  <+\infty.
\end{eqnarray*}
If  $|x'|\to +\infty,$ then
\begin{eqnarray*}
(-\Delta)^s g_\gamma(x')&=&C_{n,s} P.V. \int_{\mathbb{R}^n} \frac{g_\gamma(x')-g_\gamma(y')}{|x-y |^{n+2s}}dy \\
&=&C_{n,s} P.V. \int_{\mathbb{R}^n} \frac{g_\gamma(x')-g_\gamma(y')}{(| x'- y' |^2 +(x_1-y_1)^2)^{\frac{n+2s}{2}}}dy \\
&=&C_{n,s} P.V. \int_{\mathbb{R}^n} \frac{g_\gamma(x')-g_\gamma(y')}{(|x'- y'|^2 +| x'- y' |^2s^2)^{\frac{n+2s}{2}}}| x'- y'| dsdy'\\
&=&C_{n,s} P.V. \int_{\mathbb{R}^n} \frac{g_\gamma(x')-g_\gamma(y')}{| x'- y'|^{n-1+2s}} \cdot \frac{1}{(1+s^2)^{\frac{n+2s}{2}}} dsdy'\\
&=&C P.V. \int_{\mathbb{R}^{n-1}} \frac{g_\gamma(x')-g_\gamma(y')}{|x'- y' |^{n-1+2s}} dy'\\
&=&C P.V. \int_{\mathbb{R}^{n-1}} \frac{(1+| x '|^2)^{\frac{\gamma}{2}}-(1+|y '|^2)^{\frac{\gamma}{2}}}{|x'- y'|^{n-1+2s}} dy'\\
&\to& 0.
\end{eqnarray*}
Therefore,
\begin{eqnarray*}
| (-\Delta)^s g_\gamma(x') |\leq C.
\end{eqnarray*}
Similarly, we have
\begin{eqnarray}\label{vg}
\left | C_{n, s} P.V. \int_{\mathbb{R}^n}\frac{f(y_1)(g_\gamma(x')-g_\gamma(y'))}{| x- y|^{n+2s}}dy \right | \leq C.
\end{eqnarray}

Secondly, combining \eqref{vf} and \eqref{vg}, we derive
\begin{eqnarray*}
&&(-\Delta)^s h(x)\\
&=& C_{n, s} P.V. \int_{\mathbb{R}^n}\frac {f(x_1)g_\gamma(x')-f(y_1)g_\gamma(y')}{| x- y|^{n+2s}}dy\\
&=& C_{n, s} P.V. \int_{\mathbb{R}^n}\frac{f(x_1)g_\gamma(x')-f(y_1)g_\gamma(x')+f(y_1)g_\gamma(x')-f(y_1)g_\gamma(y')}{| x- y |^{n+2s}}dy\\
&=& C_{n, s} P.V. \int_{\mathbb{R}^n}\frac{(f(x_1)-f(y_1))g_\gamma(x')+f(y_1)(g_\gamma(x')-g_\gamma(y'))}{| x- y |^{n+2s}}dy\\
&=&g_\gamma(x')(-\Delta)^s f(x_1)+C_{n, s} P.V. \int_{\mathbb{R}^n}\frac{f(y_1)(g_\gamma(x')-g_\gamma(y'))}{| x- y|^{n+2s}}dy\\
&\geq&\frac{C}{l^{2s}}g_\gamma(x') -C.
\end{eqnarray*}
Therefore,
$$
\frac{(-\Delta)^s h(x)}{h(x)} \geq \frac{1}{h(x)}\left(\frac{C}{l^{2s}}g_\gamma(x') -C \right) \geq \frac{C}{l^{2s}}.
$$
This completes the proof of Lemma \ref{mlem1}.
\smallskip

Let
$$
\Sigma_\lambda:=\{x \in \mathbb{R}^n  \mid x_1<\lambda\},
$$
and
$$
w_\lambda (x, t)=u_\lambda (x,t )-u(x,t).
$$
Next we present the proof of narrow region principle (Theorem \ref{NR}).
\medskip

 \textbf{Proof of Theorem \ref{NR}.}

Let
  $$
  h(x)=\left[\left(1-\frac{(x_1-(\lambda-l))^2}{l^2}\right)_+^s+1\right](1+|x '|^2)^{\frac{\gamma}{2}}, \, 0<\gamma <2s,
  $$
and $$\bar{w}_\lambda(x, t)=\frac{e^{mt}w_\lambda(x,t)}{h(x)},$$
where $m>0$ is a constant.

By \eqref{EENR} and a direct calculation, $\bar{w}_\lambda(x, t)$ satisfies
  \begin{eqnarray}\label{NENR}
  \left\{\begin{array}{ll}
 \frac{\partial \bar{w}_\lambda}{\partial t} (x, t)+\frac{1}{h(x)} C_{n,s}P.V. \int_{\mathbb{R}^n}\frac{ \bar{w}_\lambda(x, t)-\bar{w}_\lambda(y, t)}{| x-y|^{n+2s}}h(y)dy\\
 =\left(c_\lambda(x, t) +m- \frac{(-\Delta)^s h(x)}{h(x)}\right)\bar{w}_\lambda(x, t), & (x, t) \in  \Omega \times \mathbb{R}, \\
 \bar{w}_\lambda(x,t)\geq 0,& (x,t )\in  (\Sigma_\lambda \backslash \Omega ) \times \mathbb{R},\\
 \bar{w}_\lambda(x,t)\to 0,& (x,t )\in \Sigma_\lambda \times \mathbb{R}, | x' |\to +\infty.
\end{array} \right.
\end{eqnarray}

Next we prove that
\begin{eqnarray}\label{barw}
\bar w_\lambda (x,t)\geq \min \left\{ 0, \inf_{\Omega}\bar{w}_\lambda(x, \bar{t})\right\},\,  \, (x,t) \in \Omega \times [\bar{t}, T],\, \forall \,[\bar{t}, T] \subset \mathbb{R}.
\end{eqnarray}
If \eqref{barw} is false, then there exists $(x_0, t_0) \in \Omega  \times (\bar{t}, T]$ such that
$$
\bar{w}_\lambda (x_0, t_0) =\inf_{\Sigma_\lambda \times (\bar{t}, T]}\bar{w}_\lambda (x, t) <\min \left\{ 0, \inf_{\Omega}\bar{w}_\lambda(x, \bar{t})\right\},
$$
hence
$$
\frac{\partial \bar{w}_\lambda}{\partial t } (x_0, t_0) \leq 0.
$$

While
\begin{eqnarray*}
&&P.V. \int_{\mathbb{R}^n}\frac{ \bar{w}_\lambda(x_0, t_0)-\bar{w}_\lambda(y, t_0)}{| x_0-y |^{n+2s}}h(y)dy\\
&=&P.V. \int_{\Sigma_\lambda}\frac{ \bar{w}_\lambda(x_0, t_0)-\bar{w}_\lambda(y, t_0)}{| x_0-y|^{n+2s}}h(y)dy+P.V. \int_{\Sigma^c_\lambda}\frac{ \bar{w}_\lambda(x_0, t_0)-\bar{w}_\lambda(y, t_0)}{| x_0-y |^{n+2s}}h(y)dy\\
&=&P.V. \int_{\Sigma_\lambda}\frac{ \bar{w}_\lambda(x_0, t_0)-\bar{w}_\lambda(y, t_0)}{| x_0-y |^{n+2s}}h(y)dy+P.V. \int_{\Sigma_\lambda}\frac{ \bar{w}_\lambda(x_0, t_0)+\bar{w}_\lambda(y, t_0)}{| x_0-y^\lambda |^{n+2s}}h(y^\lambda)dy\\
&\leq& P.V. \int_{\Sigma_\lambda}\frac{ \bar{w}_\lambda(x_0, t_0)-\bar{w}_\lambda(y, t_0)}{| x_0-y^\lambda |^{n+2s}}h(y^\lambda)dy+P.V. \int_{\Sigma_\lambda}\frac{ \bar{w}_\lambda(x_0, t_0)+\bar{w}_\lambda(y, t_0)}{|x_0-y^\lambda |^{n+2s}}h(y^\lambda)dy\\
&=& P.V. \int_{\Sigma_\lambda}\frac{ 2\bar{w}_\lambda(x_0, t_0)}{| x_0-y^\lambda |^{n+2s}}h(y^\lambda)dy\\
&<&0,
\end{eqnarray*}
where we have used the fact that
$$
| x_0-y| \leq | x_0-y^\lambda| \,\, \mbox{and}\,\,  h(y)\geq h(y^\lambda),\, y\in \Sigma_\lambda.
$$

In addition, since $c_\lambda$ is bounded from above, we apply
 \eqref{eigen} and  choose $m=\frac{C_1}{2l^{2s}},$
$$
\left(c_\lambda(x_0, t_0) +m- \frac{(-\Delta)^s h(x_0)}{h(x_0)}\right)\bar{w}_\lambda(x_0, t_0) >0,\,\, \mbox{if} \,\, l \,\,\mbox{is sufficiently small}.
$$
Therefore, we derive a contradiction from \eqref{NENR} and conclude \eqref{barw}. That is,
\begin{eqnarray*}
\bar w_\lambda (x,t)&\geq& \min \left\{ 0, \inf_{\Omega}\bar{w}_\lambda(x, \bar{t})\right\}\\
&=&\min \left\{ 0, e^{m \bar{t}}\inf_{\Omega}\frac{w_\lambda(x, \bar{t})}{h(x)}\right\}\\
&\geq& -C e^{m \bar{t}}.
\end{eqnarray*}
Therefore,
$$
w_\lambda(x, t) \geq -Ce^{-m(t-\bar{t})}h(x).
$$
Since the above inequality holds for any $\bar{t}\,(<T) \in \mathbb{R},$  let $\bar{t} \to -\infty,$ we have
$$
w_\lambda(x, t) \geq 0, \,\, (x,t) \in \Omega \times \mathbb{R}.
$$
It follows from $w_\lambda(x, t) \geq 0 $ in $(\Sigma_\lambda \backslash \Omega) \times \mathbb{R}$ that
$$
w_\lambda(x, t) \geq 0, \,\, (x,t) \in \Sigma_\lambda \times \mathbb{R}.
$$

Moreover, if there exists a point $(x^0, t^0)  \in \Omega \times \mathbb{R}$ such that
$$
w_\lambda (x^0, t^0) =0,
$$
then
$$
w_\lambda (x^0, t^0) = \inf_{\Sigma_\lambda \times \mathbb{R}} w_\lambda (x,t)=0,
$$
and
$$
\frac{\partial w_\lambda}{\partial t}  (x^0, t^0) =0, \,\, (-\Delta)^s w_\lambda(x^0, t^0)<0.
$$
Therefore, by \eqref{EENR}, we have
$$
0>\frac{\partial w_\lambda}{\partial t}  (x^0, t^0)+(-\Delta)^s w_\lambda(x^0, t^0)=c_\lambda(x^0, t^0)w_\lambda(x^0, t^0)=0,
$$
this is a contradiction and thus the strong maximum principle holds.

This completes the proof of Theorem \ref{NR}.
\medskip

Next we present the proof of maximum principle for anti-symmetric functions.

\textbf{Proof of Theorem \ref{MPA}.}

Denote
$$
\hat{a}=\sup\{|x_1|,\,\,  x\in \Omega\}.
$$
Let $x=(x_1, x') $, $a=\frac {1}{\hat{a}+1}$ and
  $$
  h(x)=\left[\left(1-a^2x_1^2\right)_+^s+1\right](1+\mid b x '\mid^2)^{\frac{\gamma}{2}},  \, \, 0<\gamma <2s.
  $$
It can be seen from the proof of Lemma \ref{mlem1} that we can choose a constant $b$ associated with $a$ in the above equality  such that
 \begin{eqnarray}\label{newh}
\frac{(-\Delta)^s h(x)}{h(x)} \geq C_1 a^{2s}.
 \end{eqnarray}

Denote
$$\bar{w}_\lambda(x, t)=\frac{e^{mt}w_\lambda(x,t)}{h(x)},$$
where $m>0$ is a constant which is determined to be later.

By \eqref{EMR} and a direct calculation, $\bar{w}_\lambda(x, t)$ satisfies
  \begin{eqnarray}\label{NEMR}
  \left\{\begin{array}{ll}
 \frac{\partial \bar{w}_\lambda}{\partial t} (x, t)+\frac{1}{h(x)} C_{n,s}P.V. \int_{\mathbb{R}^n}\frac{ \bar{w}_\lambda(x, t)-\bar{w}_\lambda(y, t)}{\mid x-y \mid^{n+2s}}h(y)dy\\
 =\left(c_\lambda(x, t) +m- \frac{(-\Delta)^s h(x)}{h(x)}\right)\bar{w}_\lambda(x, t), & (x, t) \in  \Omega \times \mathbb{R}, \\
 \bar{w}_\lambda(x,t)\geq 0,& (x,t )\in  (\Sigma_\lambda \backslash \Omega ) \times \mathbb{R},\\
 \bar{w}_\lambda(x,t)\to 0,& (x,t )\in \Sigma_\lambda \times \mathbb{R}, | x'| \to +\infty.
\end{array} \right.
\end{eqnarray}

Next we prove that
\begin{eqnarray}\label{Barw}
\bar w_\lambda (x,t)\geq \min \left\{ 0, \inf_{\Omega}\bar{w}_\lambda(x, \bar{t})\right\},\,  \, (x,t) \in \Omega \times [\bar{t}, T],\, \forall \,[\bar{t}, T] \subset \mathbb{R}.
\end{eqnarray}
If \eqref{Barw} is not valid, there exists $(x_0, t_0) \in \Omega  \times (\bar{t}, T]$ such that
$$
\bar{w}_\lambda (x_0, t_0) =\min_{\Sigma_\lambda \times (\bar{t}, T]}\bar{w}_\lambda (x, t) <\min \left\{ 0, \inf_{\Omega}\bar{w}_\lambda(x, \bar{t})\right\},
$$
then
$$
\frac{\partial \bar{w}_\lambda}{\partial t } (x_0, t_0) \leq 0.
$$

By \eqref{Clambda}, we may assume  that
$$
c_\lambda(x_0, t_0)<\frac{C_1 a^{2s}}{2},
$$
 and choose $m=\frac{C_1 a^{2s}}{2},$  we derive from \eqref{newh} and $\bar{w}_\lambda(x_0, t_0)<0$ that
$$
\left(c_\lambda(x_0, t_0) +m- \frac{(-\Delta)^s h(x_0)}{h(x_0)}\right)\bar{w}_\lambda(x_0, t_0) >0.
$$

In addition, by a similar calculation as Theorem \ref{NR}, we derive
\begin{eqnarray*}
P.V. \int_{\mathbb{R}^n}\frac{ \bar{w}_\lambda(x_0, t_0)-\bar{w}_\lambda(y, t_0)}{| x_0-y |^{n+2s}}h(y)dy
<0.
\end{eqnarray*}
Therefore, we derive a contradiction from \eqref{NEMR} and conclude \eqref{Barw}. That is,
\begin{eqnarray*}
\bar w_\lambda (x,t)&\geq& \min \left\{ 0, \inf_{\Omega}\bar{w}_\lambda(x, \bar{t})\right\}\\
&=&\min \left\{ 0, \inf_{\Omega}\frac{e^{m \bar t}w_\lambda(x, \bar{t})}{h(x)}\right\}\\
&\geq& -C e^{m \bar{t}}.
\end{eqnarray*}
It follows that
$$
w_\lambda(x, t) \geq -Ce^{-m(t-\bar{t})}h(x).
$$
Since the above inequality holds for any $\bar{t}\,(<T) \in \mathbb{R},$  let $\bar{t} \to -\infty,$ we have
$$
w_\lambda(x, t) \geq 0, \,\, (x,t) \in \Omega \times \mathbb{R},
$$
and thus conclude that
$$
w_\lambda(x, t) \geq 0, \,\, (x,t) \in \Sigma_\lambda \times \mathbb{R}.
$$

Moreover, if there exists a point $(x^0, t^0)  \in \Omega \times \mathbb{R}$ such that
$$
w_\lambda (x^0, t^0) =0,
$$
then
$$
w_\lambda (x^0, t^0) = \inf_{\Sigma_\lambda \times \mathbb{R}} w_\lambda (x,t)=0,
$$
and
$$
\frac{\partial w_\lambda}{\partial t}  (x^0, t^0) =0, \,\, (-\Delta)^s w_\lambda(x^0, t^0)<0.
$$
Therefore, by \eqref{EMR}, we have
$$
0>\frac{\partial w_\lambda}{\partial t}  (x^0, t^0)+(-\Delta)^s w_\lambda(x^0, t^0)=c_\lambda(x^0, t^0)w_\lambda(x^0, t^0)=0,
$$
this is a contradiction and thus the strong maximum principle holds.

This completes the proof of Theorem \ref{MPA}.

\section{Liouville type theorem  in  a half space}

In this section, we first prove the monotonicity of solutions to the problem
\begin{eqnarray*}\left\{\begin{array}{ll}
 \frac{\partial u}{\partial t} (x, t)+( - \Delta )^s u(x, t)=f(u(x, t)), & (x, t) \in  \mathbb{R}_+^n \times \mathbb{R}, \\
u(x,t) =0,& (x,t )\notin \mathbb{R}_+^n \times \mathbb{R}
\end{array} \right. \end{eqnarray*}
by the method of moving planes, where $f$ is a $C^1$ function satisfying $f(0)=0, \, f'(0) \leq 0.$
We also establish the connection between the existence result of the above problem and the following problem
\begin{eqnarray*}
\frac{\partial u}{\partial t} (x, t) +(-\Delta)^s u(x, t)= f(u(x, t)), \,\, (x, t) \in \mathbb{R}^{n-1} \times \mathbb{R}.
\end{eqnarray*}

Next we give the proof of Theorem \ref{MH}.
\smallskip

{\bf{Proof of Theorem \ref{MH}.}}
\medskip

(i) Let
$$
\hat{\Sigma} _\lambda =\{x \in \mathbb{R}_+^n  \mid 0< x_1< \lambda\}
$$
and
$$
w_\lambda (x, t)=u_\lambda (x,t )-u(x,t),
$$
then
\begin{eqnarray}\label{Ahe}\left\{\begin{array}{ll}
 \frac{\partial w_\lambda}{\partial t} (x, t)+( - \Delta )^s w_\lambda(x, t)=c_\lambda (x,t) w_\lambda (x,t), & (x, t) \in \hat{\Sigma}_\lambda \times \mathbb{R}, \\
w_\lambda(x,t) \geq 0,& (x,t )\in ({\Sigma} _\lambda \backslash  \hat{\Sigma} _\lambda)\times \mathbb{R},\\
w_\lambda(x^\lambda,t)=- w_\lambda(x,t),& (x,t )\in \Sigma_\lambda \times \mathbb{R}.
\end{array} \right. \end{eqnarray}
where
$$
c_\lambda (x,t) = \int_0^1 f'(su(x,t)+(1-s)u_\lambda(x,t)) ds
$$
is bounded.
\smallskip

To show the strict monotonicity of $u(x, t)$ in $x_1$-direction, we only need to prove that for any $\lambda >0$, we have
\begin{eqnarray}\label{w}
w_\lambda(x, t) > 0, \,\, (x, t) \in \Sigma_\lambda \times \mathbb{R}.
\end{eqnarray}
We divide the proof into three steps.
\medskip

{\it{Step 1.}}
We show that for the sufficiently small $\lambda$, we  have
\begin{eqnarray}\label{wgeq}
w_\lambda(x, t) \geq 0, \,\, (x, t) \in \Sigma_\lambda \times \mathbb{R}.
\end{eqnarray}
By Theorem \ref{NR} and  taking $\Omega = \hat {\Sigma}_\lambda$,  we derive
\begin{eqnarray*}
w_\lambda(x, t) \geq 0, \,\, (x, t) \in \hat{\Sigma}_\lambda \times \mathbb{R},
\end{eqnarray*}
then \eqref{wgeq} follows from the above inequality and $w_\lambda(x, t) \geq 0$ in $(\mathbb{R}^n_+)^c \times \mathbb{R}.$
\medskip

{\it{Step 2}}
\smallskip

Denote
$$
\lambda_0=\sup \{\lambda \mid w_\mu (x,t) \geq 0, \,\, (x,t )\in \Sigma_\mu \times \mathbb{R}, \,\, \mu \leq \lambda\}.
$$
In this step, we prove that $\lambda_0= +\infty.$

If not, then $0<\lambda_0<+\infty$ and there exists a sequence $\lambda_k \geq \lambda_0$ such that $\lambda_k \to \lambda_0$ and
$$
Z_k:=\{(x, t) \in \Sigma_{\lambda_k} \times \mathbb{R} \mid w_{\lambda_k} (x, t)<0\}
$$
is nonempty. Set
\begin{eqnarray*}
m_k&:=&\sup \{ u(y_1, x', t) \mid y_1 \in (0, \lambda_k), x' \in \mathbb{R}^{n-1},\,\, t \in \mathbb{R},\\
&& \mbox{ and there exists }\,  x_1 \in (0, \lambda_k)  \mbox{ such that } \, (x_1, x', t) \in Z_k\}.
\end{eqnarray*}

We consider the following two possibilities.

(a) $m_k \to 0,$

 (b) passing to a subsequence  we have $m_k \geq \varepsilon_0$  for some $\varepsilon_0 >0.$

First, assume that Case (a) holds.

It can be seen from \eqref{he} that
\begin{eqnarray}\label{lambdak}\left\{\begin{array}{ll}
 \frac{\partial w_{\lambda_k}}{\partial t} (x, t)+( - \Delta )^s w_{\lambda_k}(x, t)=c_{\lambda_k} (x,t)w_{\lambda_k}(x, t), & (x, t) \in  \hat{\Sigma}_{\lambda_k} \times \mathbb{R}, \\
w_{\lambda_k}(x,t) \geq 0,& (x,t )\in ({\Sigma} _{\lambda_k} \backslash  \hat{\Sigma} _{\lambda_k})\times \mathbb{R},\\
w_{\lambda_k}(x^{\lambda_k},t)=- w_{\lambda_k}(x,t),& (x,t )\in \Sigma_{\lambda_k} \times \mathbb{R},
\end{array} \right. \end{eqnarray}
where
$$
c_{\lambda_k} (x,t) = \int_0^1 f'(su(x,t)+(1-s)u_{\lambda_k}(x,t)) ds.
$$
Denote
$$
q_k:=\sup_{(x, t) \in Z_k} c_{\lambda_k}(x, t)
$$
By the definition of $c_{\lambda_k},$ $q_k$ and $m_k,$ and $f'(0)\leq 0$, we have
$$
\mathop {\overline{\lim} }\limits_{k\to +\infty}  q_k \leq 0.
$$
By the maximum principle for antisymmetric functions (Theorem \ref{MPA}), we have
$$
w_{\lambda_k}(x, t) \geq 0 \, \mbox{ on } Z_k.
$$
This of course contradicts the definition of $Z_k,$ and therefore Case (a) cannot occur.

Secondly, assume that Case (b) holds, then there exist subsequences $x_1^k,\, y_1^k \in (0, \lambda_k),\, z^k \in \mathbb{R}^{n-1},\, t^k \in \mathbb{R}$
such that
\begin{eqnarray}\label{wlambdak}
w_{\lambda_k} (x_1^k, z^k, t^k) <0 \,\, \mbox{and}\,\, u(y_1^k, z^k, t^k) \geq \varepsilon_0.
\end{eqnarray}
We assume that
\begin{eqnarray}\label{y1k}
x_1^k \to a,\,\, y_1^k \to b,\,\, \mbox{ for some } a, \, b \in [0, \lambda_0].
\end{eqnarray}
Consider the functions
$$
u^k(x,t):= u(x_1, x'+z^k, t+t^k), \, x=(x_1, x' ) \in \mathbb{R}^n, \, t \in \mathbb{R},
$$
and define
$$
w^k_{\lambda_k}(x, t):= u^k_{\lambda_k}(x, t) -u^k(x, t).
$$
By \eqref{wlambdak}, we have
\begin{eqnarray}\label{wklambdak}
w^k_{\lambda_k}(x_1^k, 0, 0) <0 \,\, \mbox{ and } \,\, u^k (y_1^k, 0, 0) \geq \varepsilon_0.
\end{eqnarray}

Since $u^k(x, t)$ is uniformly  bounded, by the regularity estimates for the fractional parabolic equations (\cite{FR}),  up to a subsequence (still denoted by $u^k$), as $k \to +\infty,$ we have
$$
u^k(x, t) \to \tilde u(x, t), \, \, (-\Delta)^s u^k(x, t) \to (-\Delta)^s \tilde u(x, t),
$$
and
$$
w^k_{\lambda_k}(x,t) \to  \tilde w_{\lambda_0} (x, t)=\tilde u_{\lambda_0}(x, t)-\tilde u(x, t) \geq  0, \,\, (x, t)\in \Sigma_{\lambda_0} \times \mathbb{R}.
$$
It follows that $\tilde u (x, t)$ satisfies
\begin{eqnarray}\label{tildeu}\left\{\begin{array}{ll}
 \frac{\partial \tilde u }{\partial t} (x, t)+( - \Delta )^s \tilde u (x, t)=f(\tilde u(x, t)), & (x, t) \in  \mathbb{R}^n_+ \times \mathbb{R}, \\
\tilde u (x,t) =0,& (x,t )\notin \mathbb{R}_+^n \times \mathbb{R},\\
\tilde u (x, t) \geq 0, & (x,t )\in \mathbb{R}^n \times \mathbb{R}.
\end{array} \right. \end{eqnarray}
By \eqref{y1k} and \eqref{wklambdak}, we derive that
\begin{eqnarray}\label{op}
\tilde  u ( b, 0, 0) \geq \varepsilon_0 >0.
\end{eqnarray}
Combining \eqref{tildeu} and \eqref{op}, we have
\begin{eqnarray}\label{wl}
\tilde  u ( x, t) >0, \,\, (x, t) \in \mathbb{R}_+^n \times \mathbb{R}.
\end{eqnarray}
In fact, if \eqref{wl} is false, then there exists a point $(x_0, t_0) \in \mathbb{R}_+^n \times \mathbb{R}$ such that
$$
\tilde u (x_0, t_0) =0= \inf_{\mathbb{R}^n \times \mathbb{R}}\tilde u (x, t),
$$
then
$$
\frac{\partial \tilde u}{ \partial t} (x_0, t_0) + (-\Delta)^s \tilde u(x_0, t_0) <0,
$$
which contradicts the fact that
$$
f(\tilde u(x_0, t_0)) =f(0)=0.
$$
Therefore, \eqref{wl} holds.

To proceed with the proof, we need the following Hopf's lemma for antisymmetric functions, whose proof is similar to that for Theorem 3.1 in \cite{WC}. However, for readers' convenience, we attach it in the Appendix.

Denote
$$
\tilde{\Sigma}_\lambda = \{x\in \mathbb{R}^n \mid x_1 >\lambda \},
$$
and
$$T_{\lambda}=\{x \in \mathbb{R}^n \mid x_1= \lambda\}.
$$

\begin{lemma}\label{HLA-0} (Hopf's lemma for antisymmetric functions)
Assume that  $w_\lambda(x,t) \in ( C_{loc}^{1,1}(\Omega) \cap {\cal L}_{2s})\times C^1(\mathbb{R})$ is  bounded and  satisfies
 \begin{eqnarray*}\left\{\begin{array}{ll}
 \frac{\partial w_\lambda}{\partial t} (x, t)+( - \Delta )^s w_\lambda(x, t)=c_\lambda(x, t) w_\lambda(x, t), & (x, t) \in  \tilde{\Sigma}_\lambda \times \mathbb{R}, \\
 w_\lambda(x,t)\geq 0,& (x,t )\in  \tilde{\Sigma}_\lambda \times \mathbb{R},\\
 w_\lambda(x^\lambda,t)=- w_\lambda(x,t),& (x,t )\in \tilde{\Sigma}_\lambda \times \mathbb{R},
\end{array} \right. \end{eqnarray*}
where $c_\lambda (x,t)$ is bounded from below. If there exists a point $x \in \tilde{\Sigma}_\lambda$ such that
\begin{eqnarray*}
w_\lambda(x,t)>0, \, (x,t )\in  \tilde{\Sigma}_\lambda \times \mathbb{R}.
\end{eqnarray*}
Then
$$
\frac{\partial w_\lambda}{\partial x_1} (x^0,t_0)<0, \,\, \forall\, (x^0,t_0 )\in T_\lambda \times \mathbb{R}.
$$
\end{lemma}
\medskip

Now we continue our proof in Step 2.

First, since $\tilde  u ( x, t) \equiv 0$ for $ (x, t) \in \overline{\mathbb{R}_-^n } \times \mathbb{R},$ it follows from \eqref{wl} that there exists a point $x \in \tilde{\Sigma}_{\lambda_0}$ such that
\begin{eqnarray*}
\tilde w_{\lambda_0}(x,t)>0, \, (x,t )\in  \tilde{\Sigma}_{\lambda_0} \times \mathbb{R}.
\end{eqnarray*}
Secondly, $\tilde w_{\lambda_0}(x,t)$ satisfies
 \begin{eqnarray*}\left\{\begin{array}{ll}
 \frac{\partial \tilde w_{\lambda_0}}{\partial t} (x, t)+( - \Delta )^s \tilde w_{\lambda_0}(x, t)=\tilde c_{\lambda_0}(x, t) \tilde w_{\lambda_0}(x, t), & (x, t) \in  \tilde{\Sigma}_{\lambda_0} \times \mathbb{R}, \\
 \tilde w_{\lambda_0}(x,t)\geq 0,& (x,t )\in  \tilde{\Sigma}_{\lambda_0} \times \mathbb{R},\\
\tilde w_{\lambda_0}(x^{\lambda_0},t)=- \tilde w_{\lambda_0}(x,t),& (x,t )\in \tilde{\Sigma}_{\lambda_0} \times \mathbb{R},
\end{array} \right. \end{eqnarray*}
where
$$
\tilde c_{\lambda_0} (x,t) = \int_0^1 f'(s \tilde u(x,t)+(1-s) \tilde u_{\lambda_0} (x,t)) ds.
$$
Therefore, we derive from  Hopf's lemma (Lemma \ref{HLA-0}) that
\begin{eqnarray*}
\frac{\partial\tilde w_{\lambda_0}}{\partial x_1} (x,t)<0, \,\, \forall\, (x,t )\in T_{\lambda_0} \times \mathbb{R}.
 \end{eqnarray*}

It follows that
\begin{eqnarray}\label{ux1}
\frac{\partial \tilde u}{\partial x_1} (\lambda_0, 0, 0) = -\frac{1}{2} \frac{\partial\tilde w_{\lambda_0}}{\partial x_1} (\lambda_0, 0, 0)>0.
\end{eqnarray}
Therefore, $\frac{\partial \tilde u}{\partial x_1} (x_1, 0, 0)$ is bounded from below by a positive constant in a neighborhood of $\lambda_0$ and it remains valid for $\frac{\partial  u^k}{\partial x_1} (x_1, 0, 0)$, i.e., there exists a positive constant $\delta >0$ such that for all sufficiently large $k,$
$$
\frac{\partial u}{\partial x_1}(x_1, z^k, t^k) = \frac{ \partial u^k} {\partial x_1}(x_1, 0, 0)>0, \,\, x_1 \in [\lambda_0 -\delta, \lambda_0+\delta].
$$
This contradicts the fact that
$$
w_{\lambda_k} (x_1^k, z^k, t^k) <0
$$
in \eqref{wlambdak}. Indeed, if $k$ is sufficiently large, then we have $2\lambda_k-x_1^k> x_1^k,$  and both belong to $[\lambda_0-\delta, \lambda_0+\delta]$.

Therefore, we conclude that $\lambda_0= +\infty.$.
\medskip

{\it{Step 3}}
\smallskip

Combining Step 1 with Step 2, we derive that for any $0< \lambda <+\infty,$
\begin{eqnarray*}
w_\lambda(x, t) \geq 0, \,\, (x, t) \in \Sigma_\lambda \times \mathbb{R}.
\end{eqnarray*}
It follows that if there exists a point $(x^0, t^0) \in  \Sigma_\lambda \times \mathbb{R} $ such that
$$
w_\lambda (x^0, t^0) =0 =\inf_{\mathbb{R}^n \times \mathbb{R}} w_\lambda (x, t),
$$
then
$$
\frac{\partial w_\lambda}{\partial t} (x^0, t^0)+( - \Delta )^s w_\lambda(x^0, t^0)<0,
$$
this contradicts
$$
f(u_\lambda(x^0, t^0))-f(u(x^0, t^0))=0.
$$
Therefore, for any $0< \lambda <+\infty,$
\begin{eqnarray*}
w_\lambda(x, t) > 0, \,\, (x, t) \in \Sigma_\lambda \times \mathbb{R}.
\end{eqnarray*}
It follows that
$$
\frac{\partial u}{\partial x_1} (x, t ) >0,\,\, (x, t) \in \mathbb{R}_+^n \times \mathbb{R}.
$$
\medskip

(ii)  Let $u$ be a positive bounded solution of \eqref{he}, for $k=1, 2, ... ,$ consider the functions
$$
u_k (x_1, x', t):= u(x_1+k, x', t),\,\, (x_1, x', t) \in (-k, +\infty) \times \mathbb{R}^{n-1} \times \mathbb{R}.
$$
Each of them solves the equation
$$
\frac{\partial u_k }{\partial t}(x,t) +(-\Delta)^s u_k (x, t) =f(u_k (x, t)), \,\, (x, t) \in (-k, +\infty) \times \mathbb{R}^{n-1} \times \mathbb{R}.
$$
Since the sequence $u_k(x, t)$ is uniformly bounded,
using parabolic estimates one shows that there exists a subsequence of $u_k(x, t)$ converges uniformly on arbitrary compact set to a bounded nonnegative solution $\tilde u(x, t)$ of \begin{eqnarray*}
\frac{\partial \tilde u}{\partial t} (x, t) +(-\Delta)^s \tilde u(x, t)= f(\tilde u(x, t)), \,\, (x, t) \in \mathbb{R}^{n-1} \times \mathbb{R}.
\end{eqnarray*}

 From the monotonicity result proved in (i), we derive that $u_k(x, t)$ is strict monotone increasing along $x_1$-direction in $(-k, +\infty),$  and thus $\tilde u(x, t)$ is  also monotone increasing along $x_1$-direction, by the boundedness of $\tilde u(x,t)$, we know that the limit
 $$
 \bar u(x', t)=\lim_{x_1\to +\infty} \tilde u(x_1, x', t), \,\, (x', t)  \in \mathbb{R}^{n-1}\times \mathbb{R}
 $$
 exists and satisfies
 \begin{eqnarray*}
\frac{\partial \bar u}{\partial t} (x, t) +(-\Delta)^s \bar u(x, t)= f(\bar u(x, t)), \,\, (x, t) \in \mathbb{R}^{n-1} \times \mathbb{R},
\end{eqnarray*}
which makes use of the fact that
\begin{eqnarray}\label{rd}
(-\Delta)^s_{\mathbb{R}^n} \bar u( x')=(-\Delta)^s_{\mathbb{R}^{n-1}} \bar u( x').
\end{eqnarray}
This can be proved by a direct calculations as shown in the Appendix. This proves (ii).

This completes the proof of Theorem \ref{MH}.
\medskip

Now based on Theorem \ref{MH}, we are able to derive the Liouville type theorem of the problem
\begin{eqnarray*}\left\{\begin{array}{ll}
 \frac{\partial u}{\partial t} (x, t)+( - \Delta )^s u(x, t)=u^p(x, t), & (x, t) \in  \mathbb{R}_+^n \times \mathbb{R}, \\
u(x,t) =0,& (x,t )\notin \mathbb{R}_+^n \times \mathbb{R},
\end{array} \right. \end{eqnarray*}
by applying the following nonexistence result.

\begin{lemma}\label{whh} (\cite{GK})
Let $0< \frac{\alpha n}{2s} \leq 1+\sigma$, assume $u \in( C_{loc}^{1,1} \cap {\cal L}_{2s})\times C^1(\mathbb{R} )$ satisfies
\begin{eqnarray}\label{whhe}
\left\{\begin{array}{ll}
\frac{\partial u}{\partial t} (x, t)+( - \Delta )^s u(x, t)=h(t) u^{1+\alpha}(x, t), & (x, t) \in  \mathbb{R}^n \times (0, T),\\
u(x, 0)=u_0(x) \geq 0, & x \in \mathbb{R}^n,
\end{array} \right.
\end{eqnarray}
where $u_0(x)$ is a nontrivial nonnegative  and continuous function on $\mathbb{R}^n,$ and $h(t)$ satisfies

$
(h_1) \,\, h\in C[0, +\infty), h \geq 0,
$

$
(h_2)\,\, c_0 t^\sigma \leq h(t) \leq c_1 t^\sigma \mbox{ for sufficiently large  } t, \mbox{ where } c_0, c_1 >0 \mbox{ and } \sigma >-1 \mbox{ are constants}.
$

Then the nonnegative solution $u(x,t)$ of \eqref{whhe} blows up for some $T_0>0$; and $u(x,t) =+\infty $ for every $t\geq T_0$ and $x \in \mathbb{R}^n.$
\end{lemma}

{\bf{Proof of Theorem \ref{NH}.}}
As a consequence of Lemma \ref{whh}, we derive that if $1 < p\leq \frac{n+2s}{n},$ then the following equation
\begin{eqnarray}\label{nnn}
 \frac{\partial u}{\partial t} (x, t)+( - \Delta )^s u(x, t)=u^p(x, t), & (x, t) \in  \mathbb{R}^n \times \mathbb{R},
\end{eqnarray}
has no nontrivial nonnegative bounded solution.

Combining \eqref{nnn} with the conclusion (ii) in Theorem \ref{MH}, we arrive at Theorem \ref{NH}.

\section{More relevant Liouville type theorems}

In this section, we employ the methods developed in the previous sections to prove Theorem \ref{wid}.

{\bf{Proof of Theorem \ref{wid}.}}

(i) For any given $\lambda \in \mathbb{R},$ set
$$
\Sigma_\lambda:= \{ x \in \mathbb{R}^n \mid x_1 <\lambda \},
$$
and
$$
w_\lambda(x, t)=u(x^\lambda, t)-u(x,t),
$$
then
\begin{eqnarray}\label{wlambda}
\left\{\begin{array}{ll}
 \frac{\partial w_\lambda}{\partial t} (x, t)+( - \Delta )^s w_\lambda(x, t)=c_\lambda (x,t) w_\lambda (x,t), & (x, t) \in  \Sigma_\lambda \times \mathbb{R}, \\
w_\lambda(x^\lambda,t) =-w_\lambda(x, t),& (x,t )\in  \Sigma_\lambda \times \mathbb{R},
\end{array} \right. \end{eqnarray}
where
$$
c_\lambda (x,t) = \int_0^1 f'(su(x,t)+(1-s)u_\lambda(x,t)) ds
$$
is bounded. Since $f' \leq 0,$ we have
$$c_\lambda (x,t) \leq 0,\,\, (x, t) \in \Sigma_\lambda \times \mathbb{R}.$$

We choose  the auxiliary function as
$$
g(x)=|x-(\lambda+1)e_1|^\sigma, \,\, \bar{w}_\lambda(x, t)=\frac{w_\lambda(x, t)}{g(x)},
$$
where $e_1=(1, 0, ..., 0),$ and $\sigma$ is a small positive number to be chosen as in the proof of Theorem 1 in \cite{CLZ}.

  Obviously, $\bar{w}_\lambda(x, t)$ and $w_\lambda(x, t)$ have the same sign and
  $$
  \lim_{|x|\to +\infty} \bar{w}_\lambda(x, t) =0.
  $$
For any  fixed
$t\in \mathbb{R},$ denote
$$
\bar{w}_\lambda(x(t), t)=\inf_{x \in \Sigma_\lambda}\bar{w}_\lambda(x, t).
$$

First, we conclude that for any  fixed
$t\in \mathbb{R},$ if
$$
\bar{w}_\lambda(x(t), t)<0,
$$
then
\begin{eqnarray}\label{keyineq}
\frac{\partial \bar {w}_\lambda}{\partial t} (x(t), {t}) \geq \frac{-C}{|x_1( t)-\lambda|^{2s}}\bar{w}_\lambda (x({t}), {t}).
\end{eqnarray}
In fact, by a similar calculation as (22) in \cite{CLZ}, we derive
$$
(-\Delta)^s w_\lambda (x(t), t) \leq \frac{C}{|x_1(t)-\lambda|^{2s}}w_\lambda (x(t), t).
$$
Combining this with \eqref{wlambda}, we derive
$$
\frac{\partial w_\lambda}{\partial t} (x(t), t) \geq \frac{-C}{|x_1(t)-\lambda|^{2s}}w_\lambda (x(t), t).
$$
\eqref{keyineq} follows from the above inequality and the definition of $\bar{w}_\lambda (x, t).$
\medskip

For any  fixed $t\in \mathbb{R},$ denote
$$
m(t):=\bar{w}_\lambda(x(t), t)=\inf_{x \in \Sigma_\lambda}\bar{w}_\lambda(x, t).
$$

To proceed with the proof, we need the following lemma.
\medskip

\begin{lemma}\label{conslem}
For any fixed $t \in \mathbb{R},$ if $m(t)\leq -m_0,$ then
\begin{eqnarray}\label{c0}
\frac{C}{|x_1(t)-\lambda|^{2s}} >c_0>0,
\end{eqnarray}
where $x(t)=(x_1(t), ... , x_n(t))$ is a minimum point of $\bar{w}_\lambda (x, t)$ in $\Sigma_\lambda.$
\end{lemma}
{\bf{Proof.}} If \eqref{c0} is not valid, then there exists a sequence of $\{t_k\},\, k=1,2, ...$ such that
\begin{eqnarray}\label{ck}
m(t_k) \leq -m_0,
\end{eqnarray}
and
$$
\frac{C}{|x_1(t_k)-\lambda|^{2s}} \to 0, \,\, k\to +\infty,
$$
therefore,
$$
|x_1(t_k) | \to +\infty,\,\, k\to +\infty,
$$
and it follows that
$$
m(t_k) =\bar {w}_\lambda (x(t_k), t_k) \to 0,\,\, k\to +\infty,
$$
this contradicts \eqref{ck} and hence completes the proof of Lemma \ref{conslem}.
\medskip

Now we continue our proof.

We want to show that
\begin{eqnarray}\label{keystep}
w_\lambda(x, t)\geq 0, \,\, (x, t)\in \Sigma_\lambda \times \mathbb{R},\,\, \forall\, \lambda \in \mathbb{R}.
\end{eqnarray}

If \eqref{keystep} is false, then there exits $t_0 \in \mathbb{R}$ such that
\begin{eqnarray}\label{keystep1}
-m_0:=m(t_0)=\bar{w}_\lambda(x(t_0), t_0)<0.
\end{eqnarray}

To derive a contradiction with \eqref{keystep1}, for any $\bar{t}<t_0,$ we construct a subsolution
$$
z(t)=-\bar{M} e^{-c_0(t-\bar{t})},
$$
where $c_0$ is as defined in \eqref{c0} and
$$
-\bar{M}= \inf_{\Sigma_\lambda \times \mathbb{R}} \bar{w}_\lambda(x, t).
$$

Next we prove that
\begin{eqnarray}\label{subsolu}
\bar{w}_\lambda(x, t) \geq z(t), \,\, (x, t) \in \overline{\Sigma_\lambda} \times [ \bar{t}, t_0].
\end{eqnarray}
Consider the function
$$
v(x, t)= \bar{w}_\lambda(x, t) - z(t), \,\, (x, t) \in \overline{\Sigma_\lambda} \times [ \bar{t}, t_0].
$$

For readers' convenience,   we show the definition domain of $v(x, t)$ in the following Figure 1.

\begin{center}
\begin{tikzpicture}[node distance = 0.3cm]
\draw [->, semithick] (-3.55,0) -- (4,0) node[right] {$x_1$};

 \draw [->, semithick] (1, 0,2) -- (1,0,-2) node[right] {$T_\lambda$};
 \fill [blue!50] (1,0,0)--( 1,2 ,0)--( 1,2,-2)
                                         --(1,0,-2)--cycle;
 \path (1,0) [very thick,fill=black]  circle(1.5pt) node at (1,-0.3) {$\bar t$};

 \path (-1,3.4) node at (-0.6,3.4) {$ \Sigma_\lambda\times [\bar t,t_0]$};
 \fill [orange] (-3.55,0,0)--( 0.26,0 ,0)--( -0.53,0,-2) --(-3.54,0,-2)--cycle;
 \fill [orange] (-3.55,0,0.03)--( 0.26,0 ,0.03)--( 1,0,2) --(-3.6,0,2)--cycle;
 \fill [blue!25] (-3.55,2,0)--( 1,2 ,0)--( 1,2,-2) --(-3.54,2,-2)--cycle;
 \fill [blue!25] (-3.55,2,-0.1)--( 1,2 ,-0.1)--( 1,2,2) --(-3.6,2,2)--cycle;
 \fill [blue!50](1,0,0)--( 1,2 ,0)--( 1,2,2)
                                         --(1,0,2)--cycle;
 \path (1,2,0) [very thick,fill=black]  circle(1pt) node at (1.2,1.8,0) {$ t_0$};
 \draw [->, semithick] (1,0) -- (1,3.3,0) node[right] {$t$};

 \path (-2.6 ,0,-1) node [black]{$\Sigma_\lambda$};
\node [below=0.8cm, align=flush center,text width=8cm] at (0,-0.4)
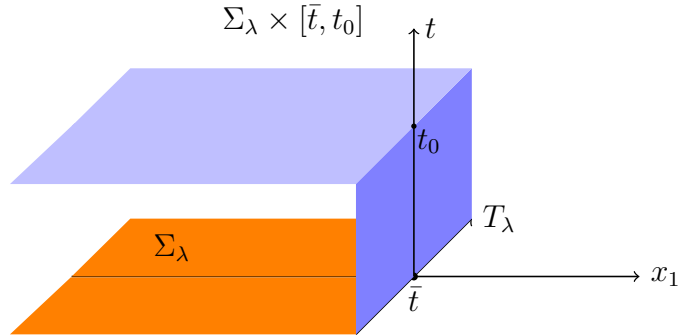
        {$Fig.$1.    \fontsize{10}{10}\selectfont  {The definition domain of $v(x, t)$.} };
\end{tikzpicture}
\end{center}

First, notice that on the bottom of the cylinder  $\Sigma_\lambda \times [ \bar{t}, t_0],$ we have
$$
v(x, t) =\bar{w}_\lambda(x, t)-z(t) =\bar{w}_\lambda(x, t)-(-\bar{M}) \geq 0,\,\, (x, t) \in \Sigma_\lambda \times \{\bar{t}\};
$$
and on the side of the domain $\Sigma_\lambda \times [ \bar{t}, t_0],$ we have
$$
v(x, t) =\bar{w}_\lambda(x, t)-z(t) =-z(t) \geq 0,\,\, (x, t) \in T_\lambda \times [ \bar{t}, t_0].
$$

Secondly, if \eqref{subsolu} does hold, then there exists a point $(x(\tilde{t}), \tilde{t}) \in \Sigma_\lambda \times (\bar{t}, t_0]$ such that
\begin{eqnarray}\label{minv}
v (x(\tilde{t}), \tilde{t}) =\inf_{\overline{\Sigma_\lambda} \times (\bar{t}, t_0]} v(x, t)<0,
\end{eqnarray}
and
\begin{eqnarray}\label{minvt}
\frac{\partial v}{\partial t} (x(\tilde{t}), \tilde{t}) \leq0,
\end{eqnarray}

On one hand, from the definition of $v(x,t),$ we have
$$
\bar{w}_\lambda(x(\tilde{t}), \tilde{t})=\inf_{{\Sigma_\lambda}} \bar{w}_\lambda(x, \tilde{t})< z (\tilde t)<0.
$$
Therefore, by \eqref{keyineq}, we have
\begin{eqnarray}\label{cc}
\frac{\partial \bar {w}_\lambda}{\partial t} (x(\tilde{t}), \tilde{t}) \geq \frac{-C}{|x_1(\tilde t)-\lambda|^{2s}}\bar{w}_\lambda (x(\tilde{t}), \tilde{t}).
\end{eqnarray}

On the other hand, we obtain from \eqref{minv} that
$$
v (x(\tilde{t}), \tilde{t}) \leq v (x(t_0), t_0),
$$
i.e.,
$$
\bar{w}_\lambda(x(\tilde{t}), \tilde{t})-\bar{w}_\lambda(x(t_0), t_0)\leq z(\tilde{t})-z(t_0)\leq 0
$$
due to the monotonicity of $z(t).$ Therefore,
\begin{eqnarray}\label{lem3condition}
m(\tilde{t})= \bar{w}_\lambda(x(\tilde{t}), \tilde{t}) \leq \bar{w}_\lambda(x(t_0), t_0)=m(t_0)=-m_0.
\end{eqnarray}
Using Lemma \ref{conslem}, we derive from \eqref{lem3condition} that
$$
\frac{C}{|x_1(\tilde{t})-\lambda|^{2s}} >c_0>0.
$$
Combining this with \eqref{cc}, we obtain
\begin{eqnarray}\label{posic}
\frac{\partial \bar {w}_\lambda}{\partial t} (x(\tilde{t}), \tilde{t}) \geq -c_0\bar{w}_\lambda (x(\tilde{t}), \tilde{t}).
\end{eqnarray}
Then by \eqref{minvt}, we derive
$$
-c_0 z(\tilde t)= \frac {\partial z} {\partial t}(\tilde t) \geq \frac {\partial \bar{w}_\lambda} {\partial t} (x(\tilde t), \tilde t) \geq -c_0 \bar{w}_\lambda (x(\tilde t), \tilde t),
$$
then
$$
z(\tilde t) \leq \bar{w}_\lambda(x(\tilde t), \tilde t),
$$
which  contradicts
$$
v (x(\tilde{t}), \tilde{t})<0.
$$
Therefore, we conclude that
$$
\bar{w}_\lambda (x, t)\geq z(t)\geq z(\bar t),\,\, (x, t) \in \overline{\Sigma_\lambda} \times [\bar{t},  t_0].
$$
i.e., we derive \eqref{subsolu}.

Let $\bar{t} \to -\infty,$ since $ z(\bar t) \to 0,$  we have
$$
\bar{w}_\lambda (x, t) \geq 0,\,\, (x, t) \in \overline{\Sigma_\lambda} \times (-\infty,  t_0].
$$
This is a contradiction with the assumption \eqref{keystep1}
$$
\bar{w}_\lambda (x(t_0), t_0) <0.
$$

As a consequence of the above results, we obtain \eqref{keystep}.

Therefore, we derive that $u(x,t)$
is  increasing in $x_1$-direction due to the arbitrariness of $\lambda,$ i.e.,
$$
\frac{\partial u}{\partial x_1} (x, t) \geq 0.
$$
Replacing $x_1$ by $-x_1,$ we have
$$
\frac{\partial u}{\partial x_1} (x, t) \leq 0.
$$
It follows that
$$
\frac{\partial u}{\partial x_1} (x, t) = 0
$$
and $u$ is a  constant for $x_1.$ By a same argument for $x_i,$ one shows that $u$ is a constant for all $x_i\, (1 \leq i\leq n).$ Thus $u$ is a solution of
$$
u_t =f(u(t)), \,\, \forall \, t \in \mathbb{R}.
$$
This proves (i).

(ii) and (iii) follows from the conclusion (i) and an elementary analysis, see Proposition 3.10 and Proposition 3.11 in \cite{X}.

\begin{corollary}\label{uid}
Assume  $f:[0, \infty) \to \mathbb{R}$ is a decreasing $C^1$ function,   $u \in( C_{loc}^{1,1} (\mathbb{R}^n_+)\cap C (\overline{\mathbb{R}^n_+}) \cap {\cal L}_{2s})\times C^1(\mathbb{R} ) $ satisfies
\begin{eqnarray}\label{uide}\left\{\begin{array}{ll}
 \frac{\partial u}{\partial t} (x, t)+( - \Delta )^s u(x, t)=f(u(x, t)), & (x, t) \in  \mathbb{R}_+^n \times \mathbb{R}, \\
u(x,t) =0,& (x,t )\notin \mathbb{R}_+^n \times \mathbb{R}.
\end{array} \right. \end{eqnarray}
Then

(i) Each positive bounded solution $u$ of \eqref{uide} is increasing in $x_1$-direction:
$$
\frac{\partial u}{\partial x_1} (x, t ) >0,\,\, (x, t) \in \mathbb{R}_+^n \times \mathbb{R}.
$$

(ii) If $f(c)<0$ as $c \geq 0$, then the nonnegative bounded solution of \eqref{uide} does not exist.
\end{corollary}
\medskip

{\bf{Proof.}} For any given $0<\lambda <+\infty,$ let
$$
\hat{\Sigma}_\lambda =\{x \in \mathbb{R}^n \mid 0< x_1< \lambda\}.
$$
By \eqref{uide}, we have
\begin{eqnarray*}\left\{\begin{array}{ll}
 \frac{\partial w_\lambda}{\partial t} (x, t)+( - \Delta )^s w_\lambda(x, t)=c_\lambda (x,t) w_\lambda (x,t), & (x, t) \in \hat{\Sigma}_\lambda \times \mathbb{R}, \\
w_\lambda(x,t) \geq 0,& (x,t )\in ({\Sigma} _\lambda \backslash  \hat{\Sigma} _\lambda)\times \mathbb{R},\\
w_\lambda(x^\lambda,t)=- w_\lambda(x,t),& (x,t )\in \Sigma_\lambda \times \mathbb{R},
\end{array} \right. \end{eqnarray*}
where
$$
c_\lambda (x,t) = \int_0^1 f'(su(x,t)+(1-s)u_\lambda(x,t)) ds
$$
is bounded. Since $f' \leq 0,$ we have
$$c_\lambda (x,t) \leq 0,\,\, (x, t) \in \hat{\Sigma}_\lambda \times \mathbb{R}.$$

Using Theorem \ref{MPA}, we derive
$$
w_\lambda(x, t) >0, \,\, (x, t) \in \hat{\Sigma}_\lambda \times \mathbb{R}.
$$
Since $0<\lambda <+\infty$ is arbitrary, we derive (i).

Since $u(x,t)$ is bounded in $\mathbb{R}^n \times \mathbb{R},$  by a similar argument as Theorem \ref{MH} (ii), the limit
$$
\bar u (x', t) = \lim_{x_1 \to  +\infty} u(x_1, x', t), \,\, (x', t) \in \mathbb{R}^{n-1} \times \mathbb{R}
$$
exists and satisfies
 \begin{eqnarray*}
\frac{\partial \bar u}{\partial t} (x, t) +(-\Delta)^s \bar u(x, t)= f(\bar u(x, t)), \,\, (x, t) \in \mathbb{R}^{n-1} \times \mathbb{R}.
\end{eqnarray*}
It follows from (iii)  in Theorem \ref{wid}  that the nonnegative bounded solution of \eqref{uide} does not exist.

This completes the proof of Corollary \ref{uid}.

\section{Appendix}

\begin{lemma}\label{HLA} (Hopf's lemma for antisymmetric functions)
Denote
$$
\tilde{\Sigma}_\lambda = \{x\in \mathbb{R}^n \mid x_1 >\lambda \}.
$$
Assume that  $w_\lambda(x,t) \in ( C_{loc}^{1,1}(\mathbb{R}^n) \cap {\cal L}_{2s})\times C^1(\mathbb{R})$ is  bounded and  satisfies
 \begin{eqnarray}\label{EHL}\left\{\begin{array}{ll}
 \frac{\partial w_\lambda}{\partial t} (x, t)+( - \Delta )^s w_\lambda(x, t)=c_\lambda(x, t) w_\lambda(x, t), & (x, t) \in  \tilde{\Sigma}_\lambda \times \mathbb{R}, \\
 w_\lambda(x,t)\geq 0,& (x,t )\in  \tilde{\Sigma}_\lambda \times \mathbb{R},\\
 w_\lambda(x^\lambda,t)=- w_\lambda(x,t),& (x,t )\in \tilde{\Sigma}_\lambda \times \mathbb{R},
\end{array} \right. \end{eqnarray}
where $c_\lambda (x,t)$ is bounded from below, if there exists a point $x \in \tilde{\Sigma}_\lambda$ such that
\begin{eqnarray}\label{WP}
w_\lambda(x,t)>0, \, (x,t )\in  \tilde{\Sigma}_\lambda \times \mathbb{R}.
\end{eqnarray}
Then
$$
\frac{\partial w_\lambda}{\partial x_1} (x^0,t_0)<0, \,\, \forall\, (x^0,t_0 )\in T_\lambda \times \mathbb{R}.
$$
\end{lemma}

\textbf{Proof.} Without loss of generality, we may assume that  $\lambda=0$ and $x^0=0$.
Let $$\tilde w(x,t)=e^{mt} w_\lambda(x,t),~m>0.$$
Since $c_\lambda (x,t)$ is bounded from below, we can  choose $m$ such that
    \begin{eqnarray} \label{633}
   m   + c _\lambda(x,t) \geq 0.
   \end{eqnarray}
 For fixed $t_0,$ $\tilde w(x, t)$  satisfies

\begin{equation}\label{eq:wcc66}
  \frac{\partial \tilde w}{\partial t}(x, t)+(-\Delta)^s \tilde w (x,t)
=   (m+c _\lambda(x,t) ) \tilde w (x,t) \geq 0,~(x,t) \in \tilde  \Sigma_\lambda  \times [t_0 -1, t_0 +1].
 \end{equation}

By \eqref{WP} and the continuity of $w_\lambda,$ there exists a set $D\subset\subset \tilde \Sigma_\lambda$ and a positive constant $c$ such that
\begin{eqnarray}\label{gc}
w_\lambda (x,t) >c, \,\,(x,t) \in   \tilde  \Sigma_\lambda  \times [t_0 -1, t_0 +1].
\end{eqnarray}
 Let    $D_\lambda$  be the reflection of   $D$ about the plane $T_\lambda$ for any time $t$.
 Denote $g(x)=x_1 \zeta(x) $, where
 $$
  \zeta(x)= \zeta(|x|)= \left\{
\begin{array}{ll} 1,~ & |x|<\varepsilon, \\
0,~& |x|\geq 2\varepsilon,\end{array}
\right.$$
 and
 $$0\leq \zeta(x) \leq 1,~  \zeta(x) \in C^\infty_0 (B_{2\varepsilon}(0)).$$
 Obviously, $g(x)$ is an antisymmetric function with respect to plane $T_0$, i.e.
 $$g(-x_1, x_2, \cdots x_n) = - g(x_1, x_2, \cdots x_n).$$

Now we construct a subsolution
$$\underline{w}(x,t)=\chi_{D\cup D_\lambda}(x) \tilde w(x,t)+\delta \eta(t)g(x),$$
where
 $$\chi_{D\cup D_\lambda}(x)=\left\{\begin{array}{ll}
1,& \quad x\in D\cup D_\lambda,\\
 0,&\quad x\not\in D \cup D_\lambda,
\end{array}
\right.$$
 and
 $\eta(t)\in C^\infty_0([t_0-1, t_0+1])$ satisfies
\begin{eqnarray*}
  \eta(t)=\left\{\begin{array}{lll} 1,& t \in [t_0-\frac{1}{2}, t_0+ \frac{1}{2}],\\
 0,&  t\not \in [t_0-1, t_0+1].
  \end{array}\right.
  \end{eqnarray*}
Since $g(x)$ is a $C^\infty_0 (B_{2\varepsilon}(0))$ function, we have $(-\Delta)^s g(x) \in C^\infty(B_{2\varepsilon}(0)).$
Then we apply the mean value theorem on $(-\Delta)^s g(x)$ at  $\bar{x} =(0, x') \in T_\lambda$ and $x=(x_1, x') \in B_{2\varepsilon}(0)$
to obtain
\begin{equation}\label{eq:w4}
|(-\Delta)^sg(x)|=|(-\Delta)^sg(\bar{x}) + \nabla((-\Delta)^s g)(\xi) \cdot (x -\bar{x})|   \leq C_0 x_1,
\end{equation}
where $\xi$ lies between $\bar{x}$ and $x$, and $(-\Delta)^sg(\bar{x})=0$ due to the antisymmetry of $g(x)$.

By the definition of fractional Laplacian  and \eqref{gc}, we derive that for each fixed $t \in [t_0-1, t_0+1]$ and for any $ x  \in  B_{2\varepsilon}(0) \cap \tilde \Sigma _\lambda$, we have
\begin{align}\label{eq:cw111}
&(-\Delta)^s(\chi_{D\cup D_\lambda}\tilde{w}(x,t))\nonumber\\
&= C_{n,s} P.V. \int_{\mathbb R^n}\frac{\chi_{D\cup D_\lambda}(x)\tilde{w}(x,t)-\chi_{D\cup D_\lambda}(y)\tilde{w}(y,t)}{|x-y|^{n+2s}}\;dy\nonumber\\
&= C_{n,s} P.V. \int_{\mathbb R^n}\frac{-\chi_{D\cup D_\lambda}(y)\tilde{w}(y,t)}{|x-y|^{n+2s}}\;dy\nonumber\\
&= C_{n,s} P.V. \int_{D}\frac{-\tilde{w}(y,t)}{|x-y|^{n+2s}}\;dy+\int_{D}\frac{-\tilde{w}(y^\lambda,t)}{|x-y^\lambda|^{n+2s}}\;dy\nonumber\\
&= C_{n,s} P.V. \int_{D}\left(\frac{1}{|x-y^\lambda|^{n+2s}}-\frac{1}{|x-y|^{n+2s}}\right)\tilde{w}(y,t)\;dy\nonumber\\
&= -C_{n, s} x_1\int_{D}\frac{2(n+2s)y_1}{\zeta(y)^{n+2s+2}}\tilde{w}(y,t)\;dy\nonumber\\
&\leq -C_1 x_1,
\end{align}
where $\zeta(y)$ is some number between $|x-y|$ and $|x-y^\lambda|$, $C_1$ is a positive constant, and the second inequality from the bottom holds due to $\lambda =0$ and the application of the mean value theorem to $h(z)=z^{-\frac{n+2s}{2}}$ over $[z_1, z_2]$ with $z_2=|x-y|^2$ and $z_2=|x-y^\lambda|^2.$

 For $(x,t)\in (B_{2\varepsilon}(0) \cap \tilde \Sigma_\lambda)\times [t_0-1,t_0+1]$, by \eqref{eq:w4} and \eqref{eq:cw111}, we obtain
$$\aligned   \frac{\partial \underline{w}}{\partial t} +(-\Delta)^s \underline{w}(x,t)
=&\delta \eta'(t) g(x)+(-\Delta)^s( \chi_{D\cup D_\lambda}\tilde{w}(x,t))+ \delta \eta(t)(-\Delta)^s g(x)
 \\
\leq& \delta \eta'(t) g(x) -C_1x_1+ \delta \eta(t) C_0x_1 .
\endaligned $$

Hence, taking $\delta$ sufficiently small, we derive
\begin{equation}\label{eq:wc33}\aligned
  \frac{\partial \underline{w}}{\partial t} +(-\Delta)^s \underline{w}(x,t) \leq0 ,~ (x,t) \in (B_{2\varepsilon}(0) \cap \tilde\Sigma_\lambda )\times [t_0-1,t_0+1].
\endaligned \end{equation}
Set $$v(x,t)=\tilde{w}(x,t)-\underline{w}(x,t).$$
Obviously, $v(x,t)=-v(x^\lambda,t)$. From  \eqref{eq:wcc66} and \eqref{eq:wc33}, we derive  that
$v(x,t)$ satisfies
\begin{equation*}
\frac{\partial v}{\partial t}(x,t)+(-\Delta)^sv(x,t) \geq0,~(x,t) \in (B_{2\varepsilon}(0) \cap\tilde \Sigma_\lambda )\times[t_0-1,t_0+1].
\end{equation*}
Also, by the definition of $\underline{w}(x,t)$, we have
$$v(x,t)\geq 0,\ \ (x,t) \in (\tilde\Sigma_\lambda\setminus (B_{2\varepsilon}(0) \cap \tilde\Sigma_\lambda))\times[t_0-1,t_0+1]$$
and $$
v(x,t_0-1)\geq 0,\ \  x\in \tilde \Sigma_\lambda.
$$
Now,  we apply the following lemma to $v(x,t)$.

\begin{lemma}(Maximum principle for antisymmetric functions \cite{CWNH})
Let $\Omega$ be a bounded domain in $\Sigma_\lambda$. Assume that   $w_\lambda(x,t)\in (C^{1,1}_{loc}(\Omega)\cap {\mathcal L}_{2s})\times C^1([0,\infty))$ is lower semi-continuous in $x$ on $\bar{\Omega}$ and satisfies
 \begin{eqnarray*}
 \left\{
\begin{array}{ll}
 \frac{\partial w_\lambda}{\partial t}(x,t)+(-\Delta)^s w_\lambda(x,t)\geq
 c_\lambda(x,t)w_\lambda(x,t), & (x,t)\in \Omega \times (0,\infty),\\
 w_\lambda(x^\lambda,t)=-w_\lambda(x,t), &  (x,t)\in  \Sigma_\lambda  \times [0,\infty),\\
     w_\lambda(x,t)\geq 0, & (x,t)\in ( \Sigma_\lambda \backslash \Omega )\times [0,\infty),\\
 w_\lambda(x,0) \geq 0, &  x \in  \Omega .
 \end{array}
\right.\end{eqnarray*}
  If $c_\lambda(x,t) $ is bounded from above, then
$$
w_\lambda(x,t)\geq 0,~(x,t)\in \Omega\times[0,T],\ \ \forall\ T>0.
$$
\end{lemma}

As an immediate consequence of this lemma with $c_\lambda(x,t)=0$, we obtain
 $$
v(x,t)\geq0,\ \ (x,t) \in  (B_{2\varepsilon}(0) \cap \tilde \Sigma_\lambda)\times[t_0-1,t_0+1].
$$
It implies that
$$e^{mt}w_\lambda(x,t)-\delta g(x)\eta(t)\geq0, \ \,\ (x,t) \in  (B_{2\varepsilon}(0) \cap \tilde\Sigma_\lambda)\times[t_0-1,t_0+1].$$
It follows that
$$ w_\lambda(x,t)\geq e^{-mt}\delta x_1 ,\ \, \  (x,t) \in  (B_{\varepsilon}(0) \cap \tilde\Sigma_\lambda)\times\left[t_0-\frac{1}{2},t_0+\frac{1}{2}\right].$$
Since $w_\lambda (0, t_0) = 0,$ we have
$$
\frac {w_\lambda(x,t_0)-w_\lambda(0, t_0)}{x_1-0} \geq \delta e^{-m t_0} >0\ \, \  x \in  B_{\varepsilon}(0) \cap \tilde\Sigma_\lambda.$$
Therefore,
$$
\frac{\partial w_\lambda}{\partial x_1} (0, t_0) >0.
$$

This completes the proof of Lemma \ref{HLA}.

\begin{lemma}\label{lemmar}
$$
(-\Delta)^s_{\mathbb{R}^n} u( x')=(-\Delta)^s_{\mathbb{R}^{n-1}} u( x').
$$
\end{lemma}

{\bf{Proof.}}
By the definition of the fractional Laplacian, we have
\begin{eqnarray}\label{rd1}
&&(-\Delta)^s_{\mathbb{R}^n} u( x')\nonumber\\
&=&C_{n, s} P.V. \int_{\mathbb{R}^n} \frac{u(x')-u(y')}{|x-y|^{n+2s}}dy\nonumber\\
&=&C_{n, s} P.V. \int_{\mathbb{R}^n} \frac{u(x')-u(y')}{(|x'-y'|^2+|x_1-y_1|^2)^{\frac{n+2s}{2}}}dy_1dy'\nonumber\\
&=&C_{n, s} P.V. \int_{\mathbb{R}^n} \frac{u(x')-u(y')}{(|x'-y'|^2+(|x'-y'|t)^2)^{\frac{n+2s}{2}}}|x'-y'|dtdy'\nonumber\\
&=&\int_{0}^{+\infty} \frac{1}{(1+t^2)^{\frac{n+2s}{2}}}dt \cdot C_{n, s} PV. \int_{\mathbb{R}^{n-1}} \frac{u(x')-u(y')}{|x'-y'|^{{n+2s-1}}}dy',
\end{eqnarray}
where
\begin{eqnarray}\label{rd2}
C_{n, s}=\left( \int_{\mathbb{R}^n} \frac{1-cos (2 \pi \zeta_1)}{|\zeta|^{n+2s}} d \zeta \right)^{-1},
\end{eqnarray}
and $\zeta =(\zeta_1, \zeta'),$ please refer to \cite{CLM}.

Since
\begin{eqnarray}\label{rd3}
&&\int_{\mathbb{R}^n} \frac{1-cos (2 \pi \zeta_1)}{|\zeta|^{n+2s}} d \zeta \nonumber \\
&=&\int_{0}^{+\infty} \frac{1}{(1+t^2)^{\frac{n+2s}{2}}}dt \int_{\mathbb{R}^{n-1}} \frac{1-cos (2 \pi \zeta_1)}{|\zeta'|^{n+2s}} d \zeta'
\end{eqnarray}
Combining \eqref{rd1}-\eqref{rd3}, we derive
$$
\int_{0}^{+\infty} \frac{1}{(1+t^2)^{\frac{n+2s}{2}}}dt \cdot C_{n, s}=C_{n-1, s}.
$$
Therefore,
$$(-\Delta)^s_{\mathbb{R}^n} u( x')= C_{n-1, s} PV. \int_{\mathbb{R}^{n-1}} \frac{u(x')-u(y')}{|x'-y'|^{{n+2s-1}}}dy'=(-\Delta)^s_{\mathbb{R}^{n-1}} u( x').
$$
This completes the proof of Lemma \ref{lemmar}.

\end{document}